\documentclass[10pt,onesite,draft]{article} 
\usepackage{amsmath,amsxtra,amssymb,latexsym, amsfonts, indentfirst}
\usepackage[mathscr]{eucal}
\newfam\cyrfam

\newenvironment{proof}{%
\par\addvspace{6pt plus3pt minus2pt}%
\noindent{\bfseries\itshape\textit Proof.}\ignorespaces} {%
\if@halmos\halmos\fi
\par\addvspace{6pt plus3pt minus2pt} }
\begin{document}
\parskip 4pt
\large
\setlength{\baselineskip}{15 truept}
\setlength{\oddsidemargin} {0.5in}
\overfullrule=0mm
\def\bfh{\vhtimeb}
\date{}

\title{\bf \large  A NOTE ON FORMULAS TRANSMUTING \\ MIXED 
MULTIPLICITIES  $^1$}

\def\b{\vntime}
\author{
 Duong Quoc Viet  and Truong Thi Hong Thanh \\ 
\small Department of Mathematics, Hanoi National University of Education\\
\small 136 Xuan Thuy Street, Hanoi, Vietnam\\
\small Emails: duongquocviet@fmail.vnn.vn \;and\; thanhtth@hnue.edu.vn\\}
  \date{}
\maketitle
\centerline{
\parbox[c]{9. cm}{
\small  ABSTRACT: This paper establishes mixed multiplicity formulas concerning the relationship between mixed multiplicities of modules and mixed multiplicities of rings via rank of modules.}}

\vspace{14pt}
\centerline{\Large\bf 1. Introduction}$\\$
 It has long been known that the mixed multiplicity is an important invariant  of Algebraic Geometry and Commutative Algebra. In 1973, Risler-Teissier \cite{Te}  showed that each mixed multiplicity of $\mathfrak{m}$-primary ideals is the multiplicity of an ideal generated by a superficial sequence. For the case of arbitrary ideals, Viet \cite{Vi1} in 2000 characterized mixed multiplicities as the  Hilbert-Samuel multiplicity via (FC)-sequences. 
In past years, the positivity and the
relationship between  mixed multiplicities  and the Hilbert-Samuel multiplicity have attracted much attention (see e.g. \cite{CP, CP1, KV, KR2,KT, MV, Re, Ro, Sw, Tr2, TV, VM, VT}).
\footnotetext[1]{\begin{itemize}
\item[ ] This research was in part supported by a grant from  NAFOSTED.
\item[ ]{\it Mathematics  Subject  Classification} (2010): Primary 13H15. Secondary 13A02, 13A15, 13A30, 14C17.  
\item[ ]$ Key\; words \;  and \; phrases:$ Noetherian  ring, mixed multiplicity,  multi-graded  module.
\end{itemize}}

In a recent paper \cite{VT1}, by a new approach, the authors   gave the additivity and reduction formulas for  mixed multiplicities of multi-graded modules  and mixed multiplicities of arbitrary ideals; and they  also showed that mixed multiplicities of arbitrary ideals are additive on exact sequences.

As a continuation, this paper  gives  mixed multiplicity formulas concerning the relationship between mixed multiplicities of modules and mixed multiplicities of rings via rank of modules.

Throughout the paper, denote by $(R,\frak{n})$  an artinian local ring with  maximal ideal $\frak{n}.$  Let $S=\bigoplus_{n_1,\ldots,n_s\ge 0}S_{(n_1,\ldots,n_s)}(s > 0)$ be  a finitely generated standard $\mathbb{N}^s$-graded algebra over $R$  and let $M=\bigoplus_{n_1,\ldots,n_s\ge 0}M_{(n_1,\ldots,n_s)}$ be  a finitely generated $\mathbb{N}^s$-graded $S$-module such that $M_{(n_1,\ldots,n_s)}=S_{(n_1,\ldots,n_s)}M_{(0,\ldots,0)}$  for all $n_1,\ldots,n_s \ge 0.$ We define  $S_{++}$ to  be $\bigoplus_{\;n_1,\ldots,n_s> 0}S_{(n_1,\ldots,n_s)}.$
Denote by $\text{Proj}\; S$ the set of the homogeneous prime ideals of $S$ which do not contain $S_{++}$. 
Set $$\text{Supp}_{++}M=\Big\{P\in \text{Proj}\; S\;|\;M_P\ne 0\Big\}$$ and
   $\dim\text{Supp}_{++}M = m.$ By \cite[Theorem 4.1]{HHRT},  $\ell_R[M_{(n_1,\ldots,n_s)}]$ is a polynomial of degree $m$ for all large $n_1,\ldots,n_s.$  
  The terms of total degree $m$ in this polynomial have the form 
$$ \sum_{k_1\:+\:\cdots\:+\:k_s\;=\;m}e(M;k_1,\ldots,k_s)\dfrac{n_1^{k_1}\cdots n_s^{k_s}}{k_1!\cdots k_s!}.$$ Then $e(M;k_1,\ldots,k_s)$ is   called the {\it  mixed multiplicity of $M$ of type $(k_1,\ldots,k_s)$ } \cite{HHRT}.
In the case that $(A, \frak m)$  is  a  noetherian   local ring with  maximal ideal $\frak{m},$ $J$ is an $\frak m$-primary ideal, $I_1,\ldots, I_s$ are ideals of $A,$ $N$ is   a finitely generated  $A$-module,  it is easily seen that 
$$ F_J(J,I_1,\ldots,I_s;N) =\bigoplus_{n_0, n_1,\ldots,n_s\ge 0}\dfrac{J^{n_0}I_1^{n_1}\cdots I_s^{n_s}{N}}{J^{n_0+1}I_1^{n_1}\cdots I_s^{n_s}{N}}$$
 is a finitely generated multi-graded $F_J(J,I_1,\ldots,I_s;A)$-module. 
Then the mixed multiplicity of $F_J(J,I_1,\ldots,I_s;N)$ of type $(k_0, k_1,\ldots,k_s)$ is denoted by $$e_A\big(J^{[k_0+1]},I_1^{[k_1]},\ldots,I_s^{[k_s]};N\big)$$ and \; called the {\it
mixed \; multiplicity of  $N$ \; with respect to ideals $J,I_1,\ldots,I_s$ of type $(k_0+1, k_1,\ldots,k_s)$} (see \cite{MV, Ve}). 

Set  $\mathrm{\bf k}= k_1,\ldots,k_s;$  $\mid\mathrm{\bf k}\mid = k_1+\cdots+k_s;$ $\mathrm{\bf I}= I_1,\ldots,I_s;$ 
 $\mathrm{\bf I}^{\mathrm{\bf n}}= I_1^{n_1},\ldots,I_s^{n_s};$ and
 $\mathrm{\bf I}^{[\mathrm{\bf k}]}= I_1^{[k_1]},\ldots,I_s^{[k_s]}.$

 Let $\mathcal{A}$ be a commutative  ring,  $\mathcal{M}$    
an $\mathcal{A}$-module, and $Q$ be the total ring of fractions of $\mathcal{A}.$ Then $\mathcal{M}$ has {\it rank}  $r$ if 
$\mathcal{M}\otimes Q$ is a free $Q$-module of rank $r.$
 
 Then we first obtain the following result for mixed multiplicities of  multi-graded modules.  

\vskip 0.2cm
\noindent
{\bf Theorem 3.1.}\;{\it Let $S$  be  a finitely generated standard \;$\mathbb{N}^s$-graded algebra  over  an  artinian local ring $R$  and  let  $M$ be a finitely generated $\mathbb{N}^s$-graded $S$-module of  positive rank such that  $S_{(1,1,\ldots,1)}$ is not contained in $ \sqrt{\mathrm{Ann}_S M}$.
 Then  $$e(M;\mathrm{\bf k})= e(S;\mathrm{\bf k})\mathrm{rank}_SM.$$}

The following theorem  and its corollaries are generalizations of classical results on the Hilbert-Samuel multiplicity (see e.g. \cite [Corollary 4.7.9]{ BH} and \cite [Corollary 11.2.6]{HS}).   

\vskip 0.2cm
\noindent
{\bf Theorem 3.4.} {\it Let  $(A, \frak m)$  be  a  noetherian   local ring  with maximal ideal $\frak m$ and   residue field $k = A/\frak{m}.$ Let  $J, I_1,\ldots,I_s$ be ideals of $A$ with $J$ being $\frak m$-primary.  Let $N$ be a finitely generated  $A$-module of positive rank. Assume  that $I=I_1\cdots I_s$  is not contained in $ \sqrt{\mathrm{Ann}_A{N}}.$  
Then we have
   $$e_A(J^{[k_0+1]},\mathrm{\bf I}^{[\mathrm{\bf k}]};N)= e_A(J^{[k_0+1]},\mathrm{\bf I}^{[\mathrm{\bf k}]}; A)\mathrm{rank}_A N.$$} 

Next, we establish  formulas  concerning the relationship between mixed multiplicities of a noetherian   local ring $A$  and mixed multiplicities of  module-finite extension rings of $A$ of positive rank
 that are  generalizations of \cite [Theorem 11.2.7]{ HS} to the mixed multiplicities of ideals. These results are started by the following theorem.

\vskip 0.2cm
\noindent
{\bf Theorem 3.9.} {\it Let  $(A, \frak m)$  be  a $d$-dimensional   noetherian   local ring  with maximal ideal $\frak{m}$ and   residue field $k = A/\frak{m}.$ Let  $J, I_1,\ldots,I_s$ be ideals of $A$ with $J$ being $\frak m$-primary. Let $B$ be a module-finite extension ring of $A$ of  positive  rank. Assume  that $I=I_1\cdots I_s$  is an ideal of positive height. Denote by  $\prod$ the set of all maximal ideals $Q$ of $B$ such that $\dim B_Q = d.$  
Set $\mathrm{\bf I}B_Q = I_1B_Q,\ldots,I_sB_Q.$  Then we have
   $$e_A(J^{[k_0+1]},\mathrm{\bf I}^{[\mathrm{\bf k}]};A)= \sum_{Q \in \; \prod}\dfrac{e_{B_Q}(JB_Q^{[k_0+1]},\mathrm{\bf I}B_Q^{[\mathrm{\bf k}]}; B_Q)[B/Q:k]}{\mathrm{rank}_A B}.$$}
\;\;\;\;The above three theorems yield  interesting consequences for the cases of domains and finite extension algebras (see Corol. 3.2; Corol. 3.3; Corol. 3.6; Corol. 3.7; Corol. 3.8; Corol. 3.12; Corol. 3.13, Section 3).
\vskip 0.2cm
    
  Our approach is based on results in \cite{VT1} and ideas in proofs
 of classical results on the Hilbert-Samuel multiplicity.  
 
This paper is divided into three  sections. Section 2 is devoted to the discussion of mixed multiplicities of multi-graded modules and the multiplicity of multi-graded Rees modules, and obtains a multiplicity formula of multi-graded Rees modules 
(see Proposition 2.1, Section 2).  Section 3 proves   mixed multiplicity formulas concerning the relationship between mixed multiplicities of modules and mixed multiplicities of rings via rank of modules.   
\\

\vspace*{24pt}
\centerline{\Large\bf2. Mixed  multiplicities  of multi-graded  modules}$\\$
\noindent   This section defines  mixed multiplicities of  multi-graded modules and  mixed multiplicities of ideals with respect to modules over local rings; and some other objects that will be used in the paper.  

Set  $\mathrm{\bf k}= k_1,\ldots,k_s;$ $\mathrm{\bf k}!= k_1!\cdots k_s!;$  $\mid\mathrm{\bf k}\mid = k_1+\cdots+k_s;$  $\mathrm{\bf n}^\mathrm{\bf k}= n_1^{k_1}\cdots n_s^{k_s}.$
 Assume that $S_{(1,\ldots,1)} \nsubseteq \sqrt{\mathrm{Ann}_SM}$ and  $\dim \mathrm{Supp}_{++}M= m,$ then   
     by \cite[Theorem 4.1]{HHRT},
    $\ell_R[M_{\mathrm{\bf n}}]$ is a polynomial of degree $m$ for all large $\mathrm{\bf n}.$ The terms of total degree $m$ in this polynomial have the form 
$$ \sum_{\mid\mathrm{\bf k}\mid\;=\;m}e(M;\mathrm{\bf k})\dfrac{\mathrm{\bf n}^\mathrm{\bf k}}{\mathrm{\bf k}!}.$$ 
Then $e(M;\mathrm{\bf k})$ are  non-negative integers not all zero, called the {\it  mixed multiplicity  of $M$  of the type ${\bf k}$} \cite{HHRT}.

From  now on, denote by 
$P_M(\mathrm{\bf n})$ the Hilbert polynomial of the Hilbert function $\ell_R[M_{\mathrm{\bf n}}]$. 

 By \cite [Proposition 2.7]{VM}, $\text{Supp}_{++}M= \emptyset$ (i.e., $S_{(1,1,\ldots,1)}\subseteq \sqrt{\mathrm{Ann}_S M}$) if and only if $M_{\mathrm{\bf n}} = 0$ for all $\mathrm{\bf n} \gg 0.$ If we assign $\dim \text{Supp}_{++}M= -\infty$ to the case that $\text{Supp}_{++}M= \emptyset$ and the degree $-\infty$ to the zero polynomial then we have
$$\deg P_M(\mathrm{\bf n})= \dim\text{Supp}_{++}M. $$

 Let  $(A, \frak m)$  be  a  noetherian   local ring  with maximal ideal $\frak{m},$ residue field $k = A/\frak{m}$ and  let $N$ be a finitely generated  $A$-module. Let $ I_1,\ldots,I_s$ be ideals of $A$ 
  such that $I_1\cdots I_s$  is not contained in $ \sqrt{\mathrm{Ann}_A{N}}$ and $J$ an $\frak m$-primary ideal.
 Put $\mathrm{\bf I}= I_1,\ldots,I_s;$ 
 $\mathrm{\bf I}^{\mathrm{\bf n}}= I_1^{n_1},\ldots,I_s^{n_s};$
 $\mathrm{\bf I}^{[\mathrm{\bf k}]}= I_1^{[k_1]},\ldots,I_s^{[k_s]}.$
 Denote by $$\frak R(\mathrm{\bf I}; A) =  \bigoplus_{n_1,\ldots,n_s\ge 0}I_1^{n_1}\cdots I_s^{n_s}$$ the multi-Rees algebra of ideals $I_1,\ldots,I_s$ and by $$\;\;\;\frak R(\mathrm{\bf I}; N) =  \bigoplus_{n_1,\ldots,n_s\ge 0}I_1^{n_1}\cdots I_s^{n_s}N$$ the multi-Rees module of ideals $I_1,\ldots,I_s$ with respect to $N.$  
 Set  
$$F_J(J,\mathrm{\bf I}; A)= \bigoplus_{n_0, n_1,\ldots,n_s\ge 0}\dfrac{J^{n_0}I_1^{n_1}\cdots I_s^{n_s}}{J^{n_0+1}I_1^{n_1}\cdots I_s^{n_s}}$$
and $$\;\;\;\;\;F_J(J,\mathrm{\bf I}; N)= \bigoplus_{n_0, n_1,\ldots,n_s\ge 0}\dfrac{J^{n_0}I_1^{n_1}\cdots I_s^{n_s}{N}}{J^{n_0+1}I_1^{n_1}\cdots I_s^{n_s}{N}}.$$
Then $F_J(J,\mathrm{\bf I};A)$ is a finitely generated standard multi-graded algebra over an  artinian local ring $A/J$ and $F_J(J,\mathrm{\bf I}; N)$ is a finitely generated multi-graded $F_J(J,\mathrm{\bf I}; A)$-module. 
Set $ I = I_1\cdots I_s$ and   
$\dim \dfrac{N}{0_N: I^\infty} = q.$  Then by \cite[Proposition 3.1]{Vi1} (see \cite[Proposition 3.1]{MV}), we get        
$$\deg P_{F_J(J,\mathrm{\bf I}; N)}(n_0, \mathrm{\bf n}) = q-1.$$ Put
$$ e\big(F_J(J,\mathrm{\bf I}; N);k_0, k_1,\ldots,k_s\big)  = e\big(J^{[k_0+1]},I_1^{[k_1]},\ldots,I_s^{[k_s]};N\big)=  e\big(J^{[k_0+1]},\mathrm{\bf I}^{[\mathrm{\bf k}]};N\big)$$ with $k_0 + \mid\mathrm{\bf k}\mid = q-1.$ Then $e\big(J^{[k_0+1]},\mathrm{\bf I}^{[\mathrm{\bf k}]};N\big) $ is called the   {\it mixed multiplicity of $N$ with respect to ideals $J,\mathrm{\bf I}$ 
of type $(k_0+1,\mathrm{\bf k})$} (see \cite{MV, Ve}).
Set $$\frak R(\mathrm{\bf I}; A)_+ =  \bigoplus_{n_1 +\cdots+n_s > 0}I_1^{n_1}\cdots I_s^{n_s}; \frak J = (J,\frak R(\mathrm{\bf I};A)_+);$$ $$
 \frak R(\mathrm{\bf I}_{\widehat{i}}\,; N) = \frak R(I_1,\ldots,I_{i-1},I_{i+1},\ldots,I_s;N); \frak J_{\widehat{i}}= (J, \frak R(\mathrm{\bf I}_{\widehat{i}}\,; A)_+)$$ and $\overline{N} = \dfrac{N}{0_N: I^\infty}.$
Then by \cite[Theorem  5.2(ii) and Note 5.4]{VT1}, we get the following.   

\vskip 0.2cm
\noindent{\bf Proposition  2.1.} We have $\sum_{k_0+\mid\mathrm{\bf k}\mid =\:q-1;\; k_i =0}e\big(J^{[k_0+1]},\mathrm{\bf I}^{[\mathrm{\bf k}]}; N\big) =  e\big(\frak J_{\widehat{i}}; \frak R(\mathrm{\bf I}_{\widehat{i}}\,;\overline{N})\big).$
 \vspace*{24pt}
  
\centerline{\Large\bf3. Some formulas for  mixed multiplicities  }

\vspace*{12pt}

\noindent
In this section, we prove the mixed multiplicity formulas concerning the relationship between mixed multiplicities of modules and mixed multiplicities of rings via rank of modules. 

First, we have the following result for $\mathbb{N}^s$-graded $S$-modules.
\vskip 0.2cm
\noindent
{\bf Theorem 3.1.}\;{\it Let $S$  be  a finitely generated standard \;$\mathbb{N}^s$-graded algebra  over  an  artinian local ring $R$  and  let  $M$ be a finitely generated $\mathbb{N}^s$-graded $S$-module of  positive rank such that  $S_{(1,1,\ldots,1)}$ is not contained in $ \sqrt{\mathrm{Ann}_S M}$. Then  $$e(M;\mathrm{\bf k})= e(S;\mathrm{\bf k})\mathrm{rank}_SM.$$}
\begin{proof}\; Since $S_{(1,1,\ldots,1)}\nsubseteq \sqrt{\mathrm{Ann}_S M},$
it \; follows \;that \; $\mathrm{Supp}_{++}M \neq \emptyset.$ 
 Let $\Lambda$\; be the set of all homogeneous prime ideals $P$\; of \; $S$ such that $P \in \mathrm{Supp}_{++}M$ $\;\text{and}\; \dim \mathrm{Proj}\;(S/P) = \dim \mathrm{Supp}_{++}M .$   Then by \cite[Theorem 3.1]{VT1}, we have  $$e(M;\mathrm{\bf k})= \sum_{P \in \Lambda}\ell(M_P)e(S/P;\mathrm{\bf k}).$$
Denote by $T$ the total ring of fractions of $S.$ Since  $M$ has   positive rank,  $M\otimes T \cong T^r $ is a free $T$-module of rank $r > 0.$ Hence 
$$M_P \cong (M\otimes T)_P \cong T_P^r \cong S_P^r \neq 0$$ for any $P \in \mathrm{Ass}\; S.$ So  $\mathrm{Ass}\;S \subseteq \mathrm{Supp}M.$ Since $S_{(1,1,\ldots,1)}\nsubseteq  \sqrt{\mathrm{Ann}_S M}$, it follows that $\emptyset \ne \Lambda \subseteq \mathrm{Min}(S/\mathrm{Ann}_SM)$  by \cite[Lemma 1.1]{HHRT}. Consequently $ \Lambda \subseteq \mathrm{Ass}\; S.$ From this it follows that $M_P  \cong S_P^r$ for any $P \in \Lambda.$ Therefore $\ell(M_P) = r\ell(S_P)$ for any $P \in \Lambda.$ This fact yields          
$$e(M;\mathrm{\bf k})= \sum_{P \in \Lambda}r\ell(S_P)e(S/P;\mathrm{\bf k})= \mathrm{rank}_S M\sum_{P \in \Lambda}\ell(S_P)e(S/P;\mathrm{\bf k}) .$$
Remember that $\mathrm{Ass}\;S \subseteq \mathrm{Supp}M.$  Hence we have  $ \dim \mathrm{Proj}\;S = \dim \mathrm{Supp}_{++}M .$ So in this case,   $\Lambda$ is also the set of all homogeneous prime ideals $P$ of  $S$ such that $P \in \mathrm{Proj}\;S$ and $$\dim \mathrm{Proj}\;(S/P) = \dim\mathrm{Proj}\;S.$$  
Therefore, $\sum_{P \in \Lambda}\ell(S_P)e(S/P;\mathrm{\bf k}) = e(S;\mathrm{\bf k})$ by \cite[Theorem 3.1]{VT1}. Thus $$e(M;\mathrm{\bf k})= e(S;\mathrm{\bf k})\mathrm{rank}_SM.\;\blacksquare$$
\end{proof}

In the case that $S$  is  a finitely generated standard graded domain  over  a  field with field of fractions $K$  and   $M$ is a finitely generated graded $S$-module, we have that $K$ is the total ring of fractions of $S$ and $K\otimes M$ is a finitely generated $K$-vector space. Hence $M$ has rank. And if $S_{(1,1,\ldots,1)}$ is not contained in $ \sqrt{\mathrm{Ann}_S M}$ then 
$M$ has positive rank: $\mathrm{rank}_SM = \mathrm{dim}_K(K\otimes M).$ We obtain the following result. 
\vskip 0.4cm
\noindent
{\bf Corollary 3.2.}\;{\it Let $S$  be  a finitely generated standard \;$\mathbb{N}^s$-graded domain  over  a  field  with field of  fractions $K$  and  let  $M$ be a finitely generated $\mathbb{N}^s$-graded $S$-module such that  $S_{(1,1,\ldots,1)}$ is not contained in $ \sqrt{\mathrm{Ann}_S M}$. Then $K\otimes M$ is a finitely generated $K$-vector space and  $$e(M;\mathrm{\bf k})= e(S;\mathrm{\bf k})\mathrm{dim}_K(K\otimes M).$$}
\;\;Now, assume that $S$  is  a finitely generated standard graded domain  over  a  field  with field of fractions $K$ and $\cal S$ is a module-finite extension domain of $S$  with field of fractions $\cal K.$ Denote by $[{\cal K} : K]$ the degree of $\cal K$ over $K.$ Then 
$K \otimes {\cal S} \cong K^{[{\cal K} : K]}$. So $\mathrm{dim}_K{(K\otimes\cal S)} = [{\cal K} : K].$ Hence by Corollary 3.2, these facts yield: 
\vskip 0.2cm
\noindent
{\bf Corollary 3.3.}\;{\it Let $S$  be  a finitely generated standard \;$\mathbb{N}^s$-graded domain  over  a  field  with field of fractions $K.$   And  suppose that $\cal S$ is a module-finite extension standard $\mathbb{N}^s$-graded domain of $S$  with field of fractions $\cal K.$ Then $[{\cal K} : K]$ is finite     
and  $$e({\cal S};\mathrm{\bf k})= e(S;\mathrm{\bf k})[{\cal K} : K].$$}
\;\;In the case of  mixed multiplicities of  ideals we have the following result.  
\vskip 0.2cm
\noindent
{\bf Theorem 3.4.} {\it Let  $(A, \frak m)$  be  a  noetherian   local ring  with maximal ideal $\frak m$ and   residue field $k = A/\frak{m}.$ Let  $J, I_1,\ldots,I_s$ be ideals of $A$ with $J$ being $\frak m$-primary. Let $N$ be a finitely generated  $A$-module of positive  rank. Assume  that $I=I_1\cdots I_s$  is not contained in $ \sqrt{\mathrm{Ann}_A{N}}.$  
Then we have
   $$e_A(J^{[k_0+1]},\mathrm{\bf I}^{[\mathrm{\bf k}]};N)= e_A(J^{[k_0+1]},\mathrm{\bf I}^{[\mathrm{\bf k}]}; A)\mathrm{rank}_AN.$$}
\noindent \begin{proof}\;\;  It is natural to suppose that one can prove Theorem 3.4 by using Theorem 3.1. But in fact, this approach seems inconvenient. The following proof of Theorem 3.4 is independent of Theorem 3.1.      

Set $\overline{N}=\dfrac{N}{0_N: I^\infty}$ and
$\bar{A}=\dfrac{A}{0: I^\infty}.$  Denote by  $\Pi$ the set of all prime ideals $\frak p $ of  $A$ such that $\frak p \in \mathrm{Min}(A/\mathrm{Ann}_A\overline{N})$ and $\dim A/\frak p = \dim \overline{N}.$ By \cite[Theorem 3.2]{VT1} we get  
$$e(J^{[k_0+1]},\mathrm{\bf I}^{[\mathrm{\bf k}]};N)= \sum_{\frak p \in \Pi}\ell({N}_{\frak p})e(J^{[k_0+1]},\mathrm{\bf I}^{[\mathrm{\bf k}]};A/\frak p).$$
Denote by $D$ the total ring of fractions of $A.$ Since $N$ has positive rank $r > 0,$   $N\otimes D$ is a free $D$-module of rank $r > 0.$ Therefore  $N\otimes D \cong D^r.$ Hence for any $\frak p \in \mathrm{Ass}A,$ we have and $N_{\frak p} \cong (N\otimes D)_{\frak p}\cong D_{\frak p}^r \cong A_{\frak p}^r \neq 0.$    So $\mathrm{Ass}A \subseteq \mathrm{Supp}N.$ 
From this it follows that 
$$\mathrm{dim}\bigg[\frac{A}{0:I^\infty}\bigg]= \mathrm{dim}\bigg[\frac{A}{\mathrm{Ann}_AN:I^\infty}\bigg]. $$
Note that  $\mathrm{Ann}_A\overline{N} = \mathrm{Ann}_AN:I^\infty,$  $$\dim \overline{N} = \mathrm{dim}\Big[\frac{A}{\mathrm{Ann}_AN:I^\infty}\Big] = \mathrm{dim}\bar A.$$
By \cite[Remark 3.3]{VT1}, we have     
$$\Pi= \Big\{\frak p \in \mathrm{Min}\Big(\frac{A}{\mathrm{Ann}_AN}\Big) \mid  \; \frak p \nsupseteq I \; \mathrm{and}\; \dim A/\frak p = \dim \overline{N} \Big\}.$$ On the other hand $\mathrm{Ass}A \subseteq \mathrm{Supp}N,$   
$$\Pi= \Big\{\frak p \in \mathrm{Min}A \mid  \; \frak p \nsupseteq I \; \mathrm{and}\; \dim A/\frak p = \dim \overline{N} \Big\}.$$
Since $\dim \overline{N} = \dim \bar {A}$ and 
$$\Big\{\frak p \in \mathrm{Min}A\mid  \; \frak p \nsupseteq I \Big\} = \mathrm{Min}\bigg[\frac{A}{0:I^\infty}\bigg], $$ 
we get
\begin{align*}
\Pi&= \Big\{\frak p \in \mathrm{Min}\bigg[\frac{A}{0:I^\infty}\bigg] \mid  \;\dim A/\frak p = \dim \overline{N}\Big\}\\
&= \Big\{\frak p \in \mathrm{Min}\bar{A} \mid  \;\dim A/\frak p = \dim \bar{A}\Big\}.
\end{align*}
It is easily seen that $\Pi \subseteq \mathrm{Ass}A,$
$$N_{\frak p} \cong (N\otimes D)_{\frak p}\cong D_{\frak p}^r \cong A_{\frak p}^r$$ for any $\frak p \in \Pi.$ Hence $\ell({N}_{\frak p}) = r \ell({A}_{\frak p})$
for any $\frak p \in \Pi.$ So we have
$$e(J^{[k_0+1]},\mathrm{\bf I}^{[\mathrm{\bf k}]};N)= r\sum_{\frak p \in \Pi}\ell({A}_{\frak p})e(J^{[k_0+1]},\mathrm{\bf I}^{[\mathrm{\bf k}]};A/\frak p).$$ Now since $$\Pi= \Big\{\frak p \in \mathrm{Min}\bar{A} \mid  \;\dim A/\frak p = \dim \bar{A}\Big\},$$ it follows that
$$\sum_{\frak p \in \Pi}\ell({A}_{\frak p})e(J^{[k_0+1]},\mathrm{\bf I}^{[\mathrm{\bf k}]};A/\frak p) = e(J^{[k_0+1]},\mathrm{\bf I}^{[\mathrm{\bf k}]};A)$$ by \cite[Theorem 3.2]{VT1}. Hence we obtain
$$e_A(J^{[k_0+1]},\mathrm{\bf I}^{[\mathrm{\bf k}]};N)= e_A(J^{[k_0+1]},\mathrm{\bf I}^{[\mathrm{\bf k}]}; A)\mathrm{rank}_A N.\;\blacksquare$$
\end{proof}

\vskip 0.2cm
  \noindent
{\bf Remark 3.5.} It would be desirable to obtain Theorem 3.4 as a consequence of Theorem 3.1, that is, to prove that  $$\mathrm{rank}_{F_J(J,\bold{I};A)} F_J(J,\bold{I};N) = \mathrm{rank}_AN.$$

Note that if \;$A$ is a domain with field of fractions\; $K$ then $N$ has rank and  
$\mathrm{rank}_AN= \dim_K (K\otimes N).$ Then as an immediate consequence of Theorem 3.4, we get the following.  
\vskip 0.2cm
\noindent
{\bf Corollary 3.6.} {\it Let  $(A, \frak m)$  be  a  noetherian local domain  with maximal ideal $\frak m$ and residue field $k = A/\frak{m},$ and field of fractions $K.$  Let  $J, I_1,\ldots,I_s$ be ideals of $A$ with $J$ being $\frak m$-primary. Let $N$ be a finitely generated  $A$-module. Assume  that $I=I_1\cdots I_s$  is not contained in $ \sqrt{\mathrm{Ann}_A{N}}.$  
Then we have
   $$e_A(J^{[k_0+1]},\mathrm{\bf I}^{[\mathrm{\bf k}]};N)= e_A(J^{[k_0+1]},\mathrm{\bf I}^{[\mathrm{\bf k}]}; A)\dim_K (K\otimes N).$$}
\;\;In particular, if $A$  is  a  domain with field of fractions $K$ and $\cal A$ is a  module-finite extension domain of $A$ with field of fractions\; $\cal K$. Then $K \otimes {\cal A} \cong K^{[{\cal K} : K]}.$ Hence   $\mathrm{dim}_K({K\otimes\cal A}) = [{\cal K}: K]$. Therefore  by Corollary 3.6, we have the following result.  
\vskip 0.2cm
\noindent
{\bf Corollary 3.7.} {\it Let  $(A, \frak m)$  be  a  noetherian   local domain  with maximal ideal $\frak m$ and  residue field $k = A/\frak{m},$ and field of fractions $K.$ Let  $J, I_1,\ldots,I_s$ be ideals of $A$ with $J$ being $\frak m$-primary. 
 Let $\cal A$ be a module-finite extension domain of $A$  with field of fractions $\cal K.$ Assume  that $I=I_1\cdots I_s \neq 0.$
Then we have
   $$e_A(J^{[k_0+1]},\mathrm{\bf I}^{[\mathrm{\bf k}]}; {\cal A})= e_A(J^{[k_0+1]},\mathrm{\bf I}^{[\mathrm{\bf k}]}; A)[{\cal K}: K].$$}

Keep the conditions as in Theorem 3.4, then we have 
$\mathrm{dim}\bigg[\frac{A}{0:I^\infty}\bigg]= \mathrm{dim}\bigg[\frac{N}{0_N:I^\infty}\bigg].$
Hence  by \cite[Corollary 2.5]{VT1}  which is a generalized version of \cite[Theorem 1.4]{Ve} and \cite[Theorem 4.4]{HHRT}, $e\Big(\big(J,\mathfrak R(\mathrm{\bf I}; A)_+\big); \mathfrak R\big(\mathrm{\bf I}; \dfrac{A}{0: I^\infty}\big)\Big)= \sum_{k_0\:+\mid\mathrm{\bf k}\mid
=\;q-1}e\big(J^{[k_0+1]},\mathrm{\bf I}^{[\mathrm{\bf k}]};A\big)$
 and 
$$e\Big(\big(J,\mathfrak R(\mathrm{\bf I}; A)_+\big); \mathfrak R\big(\mathrm{\bf I}; \dfrac{N}{0_N: I^\infty}\big)\Big)= \sum_{k_0\:+\mid\mathrm{\bf k}\mid
=\;q-1}e\big(J^{[k_0+1]},\mathrm{\bf I}^{[\mathrm{\bf k}]};N\big).$$
Therefore   
\begin{align*}
e\Big(\big(J,\mathfrak R(\mathrm{\bf I}; A)_+\big); \mathfrak R\big(\mathrm{\bf I}; \dfrac{N}{0_N: I^\infty}\big)\Big)&= \sum_{k_0\:+\mid\mathrm{\bf k}\mid
=\;q-1}e\big(J^{[k_0+1]},\mathrm{\bf I}^{[\mathrm{\bf k}]};N\big)\\ 
&=\Big[\sum_{k_0\:+\mid\mathrm{\bf k}\mid
=\;q-1}e\big(J^{[k_0+1]},\mathrm{\bf I}^{[\mathrm{\bf k}]};A\big)\Big]\mathrm{rank}_AN \\
&= e\Big(\big(J,\mathfrak R(\mathrm{\bf I}; R)_+\big); \mathfrak R\big(\mathrm{\bf I}; \dfrac{A}{0: I^\infty}\big)\Big) \mathrm{rank}_AN
 \end{align*} by Theorem 3.4.
Consequently,  
$$e\Big(\big(J,\mathfrak R(\mathrm{\bf I}; R)_+\big); \mathfrak R\big(\mathrm{\bf I}; \dfrac{N}{0_N: I^\infty}\big)\Big)= e\Big(\big(J,\mathfrak R(\mathrm{\bf I}; R)_+\big); \mathfrak R\big(\mathrm{\bf I}; \dfrac{A}{0: I^\infty}\big)\Big) \mathrm{rank}_AN.$$ 
These facts yield:
\vskip 0.2cm
\noindent
{\bf Corollary 3.8.} {\it Let  $(A, \frak m)$  be  a  noetherian   local ring  with maximal ideal $\frak m$ and   residue field $k = A/\frak{m}.$ Let  $J, I_1,\ldots,I_s$ be ideals of $A$ with $J$ being $\frak m$-primary. Let $N$ be a finitely generated  $A$-module of positive rank. Assume  that $I=I_1\cdots I_s$  is not contained in $ \sqrt{\mathrm{Ann}_A{N}}.$  
Then we have
$$e\Big(\big(J,\mathfrak R(\mathrm{\bf I}; R)_+\big); \mathfrak R\big(\mathrm{\bf I}; \dfrac{N}{0_N: I^\infty}\big)\Big)= e\Big(\big(J,\mathfrak R(\mathrm{\bf I}; R)_+\big); \mathfrak R\big(\mathrm{\bf I}; \dfrac{A}{0: I^\infty}\big)\Big) \mathrm{rank}_AN.$$}

\;\;Next, we give formulas  on  the relationship between mixed multiplicities of a noetherian   local ring $A$  and mixed multiplicities of a module-finite extension ring of $A$ of positive  rank
 that are  generalizations  of \cite [Theorem 11.2.7]{ HS} to the mixed multiplicities of ideals. 
\vskip 0.2cm
\noindent
\underline{}{\bf Theorem 3.9.} {\it Let  $(A, \frak m)$  be  a $d$-dimensional   noetherian   local ring  with maximal ideal $\frak{m}$ and   residue field $k = A/\frak{m}.$ Let  $J, I_1,\ldots,I_s$ be ideals of $A$ with $J$ being $\frak m$-primary.  Let $B$ be a   module-finite extension ring of $A$ of positive rank. Assume  that $I=I_1\cdots I_s$  is an ideal of positive height.  Denote by  $\prod$ the set of all maximal ideals $Q$ of $B$ such that $\dim B_Q = d.$  
Set $\mathrm{\bf I}B_Q = I_1B_Q,\ldots,I_sB_Q.$  Then we have
   $$e_A(J^{[k_0+1]},\mathrm{\bf I}^{[\mathrm{\bf k}]};A)= \sum_{Q \in \; \prod}\dfrac{e_{B_Q}(JB_Q^{[k_0+1]},\mathrm{\bf I}B_Q^{[\mathrm{\bf k}]}; B_Q)[B/Q:k]}{\mathrm{rank}_A B}.$$}

\vskip 0.2cm
\noindent {\bf Note 3.10:} Let $F$ be a $B$-module. Then $F$ is also an $A$-module. Assume that $F$ is an $A$-module of finite length. Since $B$ is a  module-finite extension ring of $A$, $F$ is also a $B$-module of finite length.
  Assume that $\ell_B(F) = n.$  Set $\coprod = \text{Ass}_BF.$ Then there exists a composition series of $B$-module $F$ of length $n:$ 
 $$0 = F_0 \subseteq F_1 \subseteq F_2 \subseteq\cdots \subseteq F_n=F$$ where  $F_i/F_{i-1}\cong B/P_i,$ $P_i$ maximal $(1 \le i \le n).$
And in this case, $\coprod = \{P_1,\ldots,P_n\}.$     
We have $$\ell_A(F) = \sum_{i=1}^n\ell_A(F_i/F_{i-1}) = \sum_{i=1}^n\ell_A(B/P_i).$$ Hence $\ell_A(F) $ is a sum of all the  $\ell_A(B/P) $ for $P \in \coprod,$
counted as many times as $B/P$ appears as some $B/P_i.$ This number is exactly the length of $B_P$-module $F_P.$
So $$\ell_A(F) = \sum_{P \in \coprod}\ell_A(B/P)\ell_{B_P}(F_P).$$
Because  $$\ell_A(B/P) = \dim_{k}(B/P) = [B/P : k],$$  
$$\ell_A(F) = \sum_{P \in \coprod}\ell_{B_P}(F_P)[B/P : k].\eqno(1)$$

Set $\mathbf{I}B = I_1B,\ldots,I_sB.$ Then we have $\frak R(\mathbf{I}B; B)= \frak R(\mathbf{I}; B).$   Since $B$ is a module-finite extension ring of  $A$, $\frak R(\mathbf{I}B; B)$ is a module-finite extension ring of   $\frak R(\mathbf{I}; A)$. Put $\frak N = (\frak m,\frak R(\mathrm{\bf I};A)_+).$ The notion $$e_{\frak R(\mathbf{I}; A)}((J,\frak R(\mathrm{\bf I};A)_+);\frak R(\mathbf{I}B; B))$$ will mean the Hilbert-Samuel multiplicity:   $$e_{\frak R(\mathbf{I}; A)_{\frak N}}((J,\frak R(\mathrm{\bf I};A)_+)\frak R(\mathbf{I}; A)_{\frak N};\frak R(\mathbf{I}B; B)_{\frak N}).$$ 

The proof of Theorem 3.9 is based on the following lemma. 
\vskip 0.2cm
\noindent
{\bf Lemma 3.11.}\;{\it Keeping the notation of Theorem 3.9, denote by $\Omega$ the set of all maximal homogeneous ideals $\frak M $ of $\frak R(\mathbf{I}B; B)$ such that $\dim \frak R(\mathbf{I}B; B)_\frak M = \dim \frak R(\mathbf{I}; A).$ Put $\frak J = (J,\frak R(\mathrm{\bf I};A)_+);$   $\frak J_{\frak M} = \frak J\frak R(\mathbf{I}B; B)_{\frak M}.$ Then 
$$e_{\frak R(\mathbf{I}; A)}(\frak J;\frak R(\mathbf{I}B; B))= \sum_{\frak M \in \Omega} 
e_{\frak R(\mathbf{I}B; B)_{\frak M}}(\frak J_{\frak M}; \frak R(\mathbf{I}B; B)_{\frak M})\bigg[\dfrac{\frak R(\mathbf{I}B; B)}{\frak M}:k\bigg].$$}

\begin{proof}\;\; Recall that $\frak N = (\frak m,\frak R(\mathrm{\bf I};A)_+).$ 
Denote by  $\Gamma$ the set of all maximal homogeneous ideals of $\frak R(\mathbf{I}B; B).$  Note that $\frak N \subseteq \frak M$ 
for any $\frak M \in \Gamma.$  Since $B$ is a  module-finite extension ring of $A$, $\frak R(\mathbf{I}B; B)$ is a  module-finite extension ring of $\frak R(\mathbf{I}; A).$  It is clear that   
$$\dfrac{\frak R(\mathbf{I}B; B)}{\frak J^n\frak R(\mathbf{I}B; B)}$$
is an artinian $\frak R(\mathbf{I}; A)$-module, this is also an artinian $\frak R(\mathbf{I}; B)$-module. 
 Therefore 
$$\Gamma = \text{Ass}_{\frak R(\mathbf{I}B; B)}\dfrac{\frak R(\mathbf{I}B; B)}{\frak J^n\frak R(\mathbf{I}B; B)}.$$
 Consequently,  by (1) of Note 3.10,  
$$\ell_{\frak R(\mathbf{I}; A)_{\frak N}}\dfrac{\frak R(\mathbf{I}B; B)_{\frak N}}{\frak J^n\frak R(\mathbf{I}B; B)_{\frak N}} = \sum_{\frak M \in \Gamma}\ell_{\frak R(\mathbf{I}B; B)_{\frak M}}\dfrac{\frak R(\mathbf{I}B; B)_{\frak M}}{\frak J^n\frak R(\mathbf{I}B; B)_{\frak M}}\bigg[\dfrac{\frak R(\mathbf{I}B; B)}{\frak M}:k\bigg]. \eqno(2)$$
Recall that  $\frak R(\mathbf{I}B; B)$ is a  module-finite extension ring of $\frak R(\mathbf{I}; A),$ $$\dim \frak R(\mathbf{I}B; B) = \dim \frak R(\mathbf{I}; A).$$ On the other hand ht$I > 0$,  $\dim \frak R(\mathbf{I}; A) = \dim A +s = d+s.$ 
It is easily seen that   
$\mathrm{Ann}_AB = 0,$  $\dim_A B = \dim A = d.$  Since $\mathrm{Ann}_AB = 0$ and ht$I > 0,$   
 $$\dim_{\frak R(\mathbf{I}; A)}\frak R(\mathbf{I}B; B) = \dim_AB +s = d+s.$$
So  
$$\dim \frak R(\mathbf{I}; A) = \dim_{\frak R(\mathbf{I}; A)}\frak R(\mathbf{I}B; B) = \dim  \frak R(\mathbf{I}B; B) = d+s.\eqno(3)$$
Therefore  
$$e_{\frak R(\mathbf{I}; A)}(\frak J;\frak R(\mathbf{I}B; B))=
\text{lim}_{n\rightarrow \infty}\dfrac{(d+s)!}{n^{d+s}}\ell_{\frak R(\mathbf{I}; A)_{\frak N}}\dfrac{\frak R(\mathbf{I}B; B)_{\frak N}}{\frak J^n\frak R(\mathbf{I}B; B)_{\frak N}}.\eqno(4)$$
Because that  if $\frak M \in \Gamma \setminus \Omega$ then $\dim \frak R(\mathbf{I}B; B)_{\frak M} \neq \dim \frak R(\mathbf{I}; A),$ 
by (3) it follows that 
$$\dim \frak R(\mathbf{I}B; B)_{\frak M} < d+s.\eqno(5)$$
By (2); (4) and (5) we have  
\begin{align*}
 e_{\frak R(\mathbf{I}; A)}(\frak J;\frak R(\mathbf{I}B; B))&= \text{lim}_{n\rightarrow \infty}\dfrac{(d+s)!}{n^{d+s}}\sum_{\frak M \in \Gamma}\ell_{\frak R(\mathbf{I}B; B)_{\frak M}}\dfrac{\frak R(\mathbf{I}B; B)_{\frak M}}{\frak J^n\frak R(\mathbf{I}B; B)_{\frak M}}\bigg[\dfrac{\frak R(\mathbf{I}B; B)}{\frak M}:k\bigg]\\ 
&= \text{lim}_{n\rightarrow \infty}\dfrac{(d+s)!}{n^{d+s}}\sum_{\frak M \in \Omega}\ell_{\frak R(\mathbf{I}B; B)_{\frak M}}\dfrac{\frak R(\mathbf{I}B; B)_{\frak M}}{\frak J^n\frak R(\mathbf{I}B; B)_{\frak M}}\bigg[\dfrac{\frak R(\mathbf{I}B; B)}{\frak M}:k\bigg]. \end{align*} Thus, 
$$e_{\frak R(\mathbf{I}; A)}(\frak J;\frak R(\mathbf{I}B; B))= \sum_{\frak M \in \Omega} 
e_{\frak R(\mathbf{I}B; B)_{\frak M}}(\frak J_{\frak M}; \frak R(\mathbf{I}B; B)_{\frak M})\bigg[\dfrac{\frak R(\mathbf{I}B; B)}{\frak M}:k\bigg].\;\blacksquare$$
\vskip 0.2cm
\noindent
{\bf The proof of Theorem 3.9:}
It is easily seen that there is an one-to-one   correspondence  
between the set of maximal ideals $\Omega$  and the set of maximal ideals  $\prod$ given by
$$\frak M \mapsto 
Q = \frak M \bigcap B.$$
 Moreover, if $Q = \frak M \bigcap B$ then $\dfrac{\frak R(\mathbf{I}B; B)}{\frak M}\cong  \dfrac{B }{Q}$  is a finite extension field of $k$ and
$$\frak R(\mathbf{I}B; B)_{\frak M}= \frak R(\mathbf{I}B_Q; B_Q)_{\frak M} \;\text{and}\; \frak J_{\frak M} = (JB_Q,\frak R(\mathrm{\bf I}B_Q;B_Q)_+)\frak R(\mathbf{I}B_Q; B_Q)_{\frak M}. \eqno(6)$$
 Remember \; that \; $\mathrm{\bf I}^{\mathrm{\bf u}}= I_1^{u_1},\ldots,I_s^{u_s}$\;\; and \;\; 
 $\mathrm{\bf I}^{\mathrm{\bf u}}B_Q= I_1^{u_1}B_Q,\ldots,I_s^{u_s}B_Q$
\;\;for \;\;any \;\;$\mathrm{\bf u}= u_1,\ldots,u_s.$ 
 By \cite[Proposition 2.4]{VT1}  which is a generalized version of \cite[Theorem 1.4]{Ve} and \cite[Theorem 4.4]{HHRT}, we have
$$e_{\frak R(\mathrm{\bf I}^\mathrm{\bf u}; A)}\big(\big(J,\frak R(\mathrm{\bf I}^\mathrm{\bf u}; A)_+\big); \frak R(\mathrm{\bf I}^\mathrm{\bf u}; B)\big)= \sum_{k_0\:+\mid\mathrm{\bf k}\mid\;
=\;d-1}e_A\big(J^{[k_0+1]},\mathrm{\bf I}^{[\mathrm{\bf k}]}; B\big)\mathrm{\bf u}^\mathrm{\bf k}$$ and
\begin{align*}
&e_{\frak R(\mathrm{\bf I}^\mathrm{\bf u}B_Q; B_Q)}\big(\big(JB_Q,\frak R(\mathrm{\bf I}^\mathrm{\bf u}B_Q; B_Q)_+\big); \frak R(\mathrm{\bf I}^\mathrm{\bf u}B_Q; B_Q)\big)\\
&= \sum_{k_0\:+\mid\mathrm{\bf k}\mid\;
=\;d-1}e_{B_Q}\big(JB_Q^{[k_0+1]},\mathrm{\bf I}B_Q^{[\mathrm{\bf k}]}; B_Q\big)\mathrm{\bf u}^\mathrm{\bf k}.\end{align*}
Hence by Lemma 3.11 and (6), we obtain
\begin{align*}
&\sum_{k_0\:+\mid\mathrm{\bf k}\mid\;
=\;d-1}e_A\big(J^{[k_0+1]},\mathrm{\bf I}^{[\mathrm{\bf k}]}; B\big)\mathrm{\bf u}^\mathrm{\bf k}\\
&= e_{\frak R(\mathrm{\bf I}^\mathrm{\bf u}; A)}\big((J,\frak R(\mathrm{\bf I}^\mathrm{\bf u}; A)_+); \frak R(\mathrm{\bf I}^\mathrm{\bf u}B; B)\big)\\
&= \sum_{Q \in \prod}e_{\frak R(\mathrm{\bf I}^\mathrm{\bf u}B_Q; B_Q)}\big(\big(JB_Q,\frak R(\mathrm{\bf I}^\mathrm{\bf u}B_Q; B_Q)_+\big); \frak R(\mathrm{\bf I}^\mathrm{\bf u}B_Q; B_Q)\big)[B/Q:k]\\
\mathrm{}&= \sum_{k_0\:+\mid\mathrm{\bf k}\mid\;
=\;d-1}\Big[\sum_{Q \in \prod} e_{B_Q}\big(JB_Q^{[k_0+1]},\mathrm{\bf I}B_Q^{[\mathrm{\bf k}]}; B_Q\big)[B/Q:k]\Big]\mathrm{\bf u}^\mathrm{\bf k}.
\end{align*} From this it follows that
$$e_A(J^{[k_0+1]},\mathrm{\bf I}^{[\mathrm{\bf k}]};B)= \sum_{Q \in \prod}e_{B_Q}(JB_Q^{[k_0+1]},\mathrm{\bf I}B_Q^{[\mathrm{\bf k}]}; B_Q)[B/Q:k]. $$
Remember that by Theorem 3.4, $$e_A(J^{[k_0+1]},\mathrm{\bf I}^{[\mathrm{\bf k}]}; B)= e_A(J^{[k_0+1]},\mathrm{\bf I}^{[\mathrm{\bf k}]}; A)\mathrm{rank}_A B.$$
Thus,
$$e_A(J^{[k_0+1]},\mathrm{\bf I}^{[\mathrm{\bf k}]};A)= \sum_{Q \in \; \prod}\dfrac{e_{B_Q}(JB_Q^{[k_0+1]},\mathrm{\bf I}B_Q^{[\mathrm{\bf k}]}; B_Q)[B/Q:k]}{\mathrm{rank}_A B}. \;\blacksquare$$

\end{proof}

When $A$ is a domain with field of fractions $K$ and $B$ is a   module-finite extension ring of $A$,  $B$ has rank and  
$\mathrm{rank}_AB = \dim_K (K \otimes B).$ Then  we obtain the following result by  Theorem 3.9. 
\vskip 0.2cm
\noindent
{\bf Corollary 3.12.} {\it Let  $(A, \frak m)$  be  a $d$-dimensional   noetherian   local domain  with maximal ideal $\frak{m}$ and  residue field $k = A/\frak{m},$ and field of fractions $K.$  Let  $J, I_1,\ldots,I_s$ be ideals of $A$ with $J$ being $\frak m$-primary.
 Let $B$ be a   module-finite extension ring of $A$.     Denote by  $\prod$ the set of all maximal ideals $Q$ of $B$ such that $\dim B_Q = d.$  
 Assume  that $I=I_1\cdots I_s \neq 0.$ Then we have
   $$e_A(J^{[k_0+1]},\mathrm{\bf I}^{[\mathrm{\bf k}]};A)= \sum_{Q \in \; \prod}\dfrac{e_{B_Q}(JB_Q^{[k_0+1]},\mathrm{\bf I}B_Q^{[\mathrm{\bf k}]}; B_Q)[B/Q:k]}{\mathrm{dim}_K (K\otimes B)}.$$}

Let $A$ be a domain with field of fractions $K,$ and $B$  a  module-finite extension domain of $A$ with field of fractions $\cal K.$ Then  $\mathrm{dim}_K (K \otimes B) = [{\cal K} : K].$ Hence  as an immediate consequence of Corollary 3.12, we give the following result.
\vskip 0.2cm
\noindent
{\bf Corollary 3.13.} {\it Let  $(A, \frak m)$  be  a $d$-dimensional   noetherian   local domain  with maximal ideal $\frak{m}$ and  residue field $k = A/\frak{m},$ and field of fractions $K.$  Let  $J, I_1,\ldots,I_s$ be ideals of $A$ with $J$ being $\frak m$-primary.
 Let $B$ be a  module-finite extension domain of $A$ with field of fractions ${\cal K}.$ Assume  that $I=I_1\cdots I_s \neq 0.$   Denote by  $\prod$ the set of all maximal ideals $Q$ of $B$ such that $\dim B_Q = d.$  
  Then we have
   $$e_A(J^{[k_0+1]},\mathrm{\bf I}^{[\mathrm{\bf k}]};A)= \sum_{Q \in \; \prod}\dfrac{e_{B_Q}(JB_Q^{[k_0+1]},\mathrm{\bf I}B_Q^{[\mathrm{\bf k}]}; B_Q)[B/Q:k]}{[{\cal K} : K]}.$$}

\vskip 0.2cm
\noindent
{\bf Acknowledgement:} Special thanks are due to the referee whose remarks substantially improved this  paper.

\end{document}